\documentclass[preprint,12pt]{elsarticle}

\usepackage{graphicx}
\usepackage{enumitem}
\usepackage{caption}
\usepackage{subcaption}
\usepackage{float}
\usepackage{booktabs}
\usepackage{siunitx}
\usepackage{microtype}
\usepackage{tikz}
\usetikzlibrary{arrows.meta,decorations.pathreplacing}

\usepackage{amsmath, amssymb}
\usepackage{bm}
\usepackage{hyperref}
\usepackage{xcolor}
\newcommand{\rev}[1]{{\color{black}#1}}
\newcommand{\revcolor}{\color{black}}

\journal{Engineering Analysis with Boundary Elements}

\begin{document}

\begin{frontmatter}

\title{Unified Regularization of 2D Singular Integrals for Axisymmetric Galerkin BEM in Eddy-Current Evaluation}

\author[whu]{Yao Luo\corref{cor1}}
\address[whu]{School of Electrical Engineering and Automation, Wuhan University, Wuhan 430072, China}

\cortext[cor1]{Corresponding author. E-mail: \texttt{sturmjungling@gmail.com}. 
Currently visiting researcher at the Department of Electronics, Information and Bioengineering (DEIB), Politecnico di Milano.}

\begin{abstract}
This paper presents an axisymmetric Galerkin boundary element method (BEM) for modeling eddy-current interactions between excitation coils and conductive objects. The formulation derives boundary integral equations from the Stratton–Chu representation for the azimuthal component of the vector potential in both air and conductive regions. The central contribution is a unified 
regularization framework for the two-dimensional (2D) singular integrals arising in Galerkin BEM. This framework handles both logarithmic and Cauchy singularities through a common set of integral transformations, eliminating the need for case-by-case analytical singularity extraction and enabling straightforward numerical quadrature. The regularization and quadrature stability are proved and verified numerically. The method is validated on several representative axisymmetric geometries, including cylindrical, conical, and spherical shells. Numerical experiments demonstrate consistently high accuracy and computational efficiency \rev{over the tested frequency interval} and coil lift-off distances. The results confirm that the proposed axisymmetric Galerkin BEM, combined with the integral transformation technique, provides a robust and efficient framework for axisymmetric eddy-current nondestructive evaluation.
\end{abstract}


\end{frontmatter}

\section{Introduction}
Eddy-current nondestructive evaluation (EC-NDE) is widely used for conductivity measurement, flaw detection, thickness estimation, and material sorting \cite{GarciaMartin2011,Bowler2019}. Many important inspection scenarios, such as heat-exchanger tubes, steam-generator tubing, cylindrical shells, and rod-type components, exhibit axial symmetry. In such cases, boundary-integral formulations are especially attractive because only the one-dimensional meridional curve needs to be discretized, whereas finite-element method (FEM) requires meshing of the entire 2D cross-section. \rev{This dimensional reduction is particularly valuable in frequency sweeps, inverse identification, and geometry-parameter studies, where the forward model must be evaluated repeatedly.}

Most existing axisymmetric BEM models employ collocation formulations \cite{Bakr1986}. Despite the conceptual simplicity of this framework, several computational challenges arise in practice. While the axisymmetric Laplace kernel admits closed-form expressions in terms of elliptic integrals, the axisymmetric Helmholtz kernel does not \cite{Tsuchimoto1990}. Consequently, several studies have developed ad hoc numerical techniques to approximate the Helmholtz kernel \cite{Cheng1992, Priede2006}. Moreover, even for the axisymmetric Laplace kernel, the closed-form expressions for singular integrals become intricate for constant elements, and their extension to higher-order elements dramatically increases the algebraic complexity \cite{Dawson1995}. \rev{While effective for the constant and linear elements, this strategy is inherently tied to the kernel-specific singular decompositions, geometry-specific formulae, explicitly derived coefficient sets, additional branching logic, and limiting configurations near singular points. For curved isoparametric elements, the affine geometric structure that enables the closed-form extraction is lost, and analytical extraction no longer provides a practical singular-integration formula. These same difficulties arise in the analytical extraction of singularities from the Helmholtz kernel, since the singular part corresponds precisely to the Laplace kernel.}

The Galerkin boundary element method (Galerkin BEM) offers a natural alternative that is particularly well-suited for higher-order accuracy, as its variational foundation provides stable projection onto the trial space and guarantees quasi-optimal error bounds \cite{Steinbach2008}. However, its application to axisymmetric eddy-current analysis has remained limited. The primary obstacle is numerical rather than conceptual: Galerkin BEM requires reliable evaluation of double integrals involving products of basis functions and singular kernels. \rev{Recent work on Galerkin BEM confirms that the treatment of singular element-pair integrals remains an active topic~\cite{Seibel2023,Montanelli2022,Montanelli2024}.}

The present work addresses this gap by introducing a coordinate transformation tailored to 2D Galerkin BEM. This construction yields smooth integrands in all relevant configurations, enabling direct evaluation via Gauss-Legendre (GL) quadrature without the analytical extraction linked to the kernel or element order. \rev{The proposed scheme regularizes logarithmic and Cauchy singularities directly on the parameter domain and works unchanged for straight and curved boundary elements. Since this treatment relies solely on the singular nature of the kernels and is independent of the physical model, it can be applied to any boundary integral equation (BIE) whose Galerkin discretization leads to 2D logarithmic or Cauchy-singular integrals.}

\rev{The remainder of this paper is organized as follows. Section~2 presents the axisymmetric Galerkin BEM formulation for eddy-current evaluation. Section~3 introduces the unified coordinate transformations for the logarithmic and Cauchy-singular integrals arising from coincident and touching element pairs. Section~4 validates the method through benchmark problems involving cylindrical, conical, and spherical conductors. Finally, Section~5 summarizes the conclusions.}

\section{Boundary integral equations for the axisymmetric vector potential}

\subsection{Boundary integral equations derived from the Stratton-Chu formula}

Consider an eddy-current coil carrying current with amplitude $I$ and angular frequency $\omega$, placed near a conductive object with conductivity $\sigma$ and permeability $\mu$. Neglecting the displacement current, the governing equations for the vector potential $\bm{A}$ are

\begin{equation}
\nabla^2 \bm{A} = -\mu_0 \bm{J}^{(e)}
\label{eq:air_vector_poisson}
\end{equation}
in the air, where $\bm{J}^{(e)}$ is the current density of the source, and

\begin{equation}
\nabla^2 \bm{A} + k^2 \bm{A} = 0
\label{eq:conductor_vector_helmholtz}
\end{equation}
in the conductive domain, where
\begin{equation}
k^2 = -i \omega \mu \sigma
\label{eq:complex_wavenumber}
\end{equation}
Here $k$ denotes the complex wavenumber with $\Re(k)>0$. BIEs will be derived with the vector version of Green’s second identity, i.e. the Stratton-Chu formula \cite{Stratton_Chu1939}:

\begin{equation}
\begin{split}
& \int_{\Omega} \Big[G(\bm{x},\bm{y}) \nabla^2 \bm{A}-\bm{A} \nabla^2 G(\bm{x},\bm{y}) \Big] d\Omega_{\bm{y}} \\
&= -\int_{\Gamma} \Big[(\nabla \cdot \bm{A})G(\bm{x},\bm{y})\bm{n} - (\bm{n} \cdot \bm{A})\nabla G(\bm{x},\bm{y}) \\
& \quad - \bm{n}\times (\nabla\times \bm{A}) G(\bm{x},\bm{y}) - (\bm{n}\times \bm{A}) \times \nabla G(\bm{x},\bm{y}) \Big] d\Gamma_{\bm{y}}
\end{split}
\label{eq:stratton_chu_laplace}
\end{equation}
for the free space, where $\Gamma$ denotes the boundary of the region $\Omega$, $\bm{n}$ is the unit normal vector pointing into the free space, and $G(\bm{x},\bm{y})$ is the fundamental solution (FS)

\begin{equation}
G(\bm{x},\bm{y}) = \frac{1}{4\pi |\bm{x}-\bm{y}|}
\label{eq:laplace_fs}
\end{equation}
satisfying

\begin{equation}
\nabla^2 G(\bm{x},\bm{y}) = -\delta(\bm{x}-\bm{y})
\label{eq:laplace_fs_identity}
\end{equation}
where $\bm{x}$ and $\bm{y}$ are the \rev{observation} and source points, respectively, $\delta$ is the Dirac delta function, and $|\bm{x}-\bm{y}|$ is the distance between them.

When $\bm{x}$ is on the boundary $\Gamma$, a BIE can be obtained from \eqref{eq:air_vector_poisson}, \eqref{eq:stratton_chu_laplace}, and \eqref{eq:laplace_fs_identity} by a limiting procedure:

\begin{align}
\frac{\Omega(\bm{x})}{4\pi} \bm{A}(\bm{x})\,\rev{-}\int_{\Gamma} \Big[(\bm{n}\cdot\bm{A})\nabla G(\bm{x},\bm{y})+ \bm{n}\times (\nabla \times \bm{A})G(\bm{x},\bm{y}) \nonumber \\
 + (\bm{n}\times\bm{A})\times\nabla G(\bm{x},\bm{y}) \Big] d\Gamma_{\bm{y}} 
= \bm{A}^{(e)}(\bm{x}),
\label{eq:bie_air_vector}
\end{align}
where $\Omega(\bm{x})$ is the solid angle at point $\bm{x}$, and $\bm{A}^{(e)}(\bm{x})$ is the source vector potential at the same point:

\begin{equation}
\bm{A}^{(e)}(\bm{x}) = \mu_0 \int_{\Omega} \bm{J}^{(e)}(\bm{y}) G(\bm{x},\bm{y}) d\Omega_{\bm{y}}.
\label{eq:source_vector_potential}
\end{equation}
In \eqref{eq:bie_air_vector}, the Coulomb gauge $\nabla \cdot \bm{A}=0$ is assumed.
Applying \eqref{eq:conductor_vector_helmholtz} together with the Helmholtz counterpart of \eqref{eq:stratton_chu_laplace} yields
\begin{align}
& \int_{\Omega} \Big[G_k(\bm{x},\bm{y})(\nabla^2\bm{A}+k^2\bm{A}) - \bm{A}(\nabla^2G_k(\bm{x},\bm{y})+k^2G_k(\bm{x},\bm{y}))\Big] d\Omega_{\bm{y}} \nonumber \\
= & \int_{\Gamma} \Big[ \bm{n} (\nabla \cdot \bm{A}) G_k(\bm{x},\bm{y}) - (\bm{n}\cdot\bm{A})\nabla G_k(\bm{x},\bm{y}) \nonumber \\
& - \bm{n}\times (\nabla\times \bm{A}) G_k(\bm{x},\bm{y}) - (\bm{n}\times\bm{A})\times \nabla G_k(\bm{x},\bm{y}) \Big] d\Gamma_{\bm{y}},
\label{eq:stratton_chu_helmholtz}
\end{align}
where $G_k(\bm{x},\bm{y}) = \exp(-ik|\bm{x}-\bm{y}|)/(4\pi|\bm{x}-\bm{y}|)$ denotes the FS satisfying
\begin{equation}
\nabla^2 G_k + k^2 G_k = -\delta(\bm{x}-\bm{y}).
\label{eq:helmholtz_fs_identity}
\end{equation}
Consequently, the BIE for the conductive region becomes
\begin{align}
\frac{\Omega(\bm{x})}{4\pi} \bm{A}(\bm{x})& +\int_{\Gamma} \Big[(\bm{n}\cdot\bm{A})\nabla G_k(\bm{x},\bm{y})+ \bm{n}\times (\nabla \times \bm{A})G_k(\bm{x},\bm{y}) \nonumber \\
 & + (\bm{n}\times\bm{A})\times\nabla G_k(\bm{x},\bm{y}) \Big] d\Gamma_{\bm{y}} 
= 0.
\label{eq:bie_conductor_vector}
\end{align}

For an axisymmetric problem, the vector potential $\bm{A}$ possesses only an azimuthal component in cylindrical coordinates $(r,\varphi,z)$. Consequently, \eqref{eq:bie_air_vector} and \eqref{eq:bie_conductor_vector} reduce to the axisymmetric BEM formulation
\begin{equation}
c(\bm{x}) u^{(1)}(\bm{x}) + \int_{\gamma} \left[\mathcal{G}(\bm{x},\bm{y})q^{(1)}(\bm{y})-u^{(1)}(\bm{y})\partial_{\bm{n}}\mathcal{G}(\bm{x},\bm{y})\right] r_{\bm{y}}ds_{\bm{y}} = u^{(e)}(\bm{x})
\label{eq:axisym_bie_air}
\end{equation}
and
\begin{equation}
\left(1-c(\bm{x})\right) u^{(2)}(\bm{x}) + \int_{\gamma} \left[u^{(2)}(\bm{y})\partial_{\bm{n}}\mathcal{G}_k(\bm{x},\bm{y})-\mathcal{G}_k(\bm{x},\bm{y})q^{(2)}(\bm{y}) \right] r_{\bm{y}}ds_{\bm{y}} = 0,
\label{eq:axisym_bie_conductor}
\end{equation}
where the superscripts (1) and (2) denote the non-conductive and conductive regions, respectively, $\gamma$ denotes the meridional curve, and
\begin{equation}
\begin{split}
u(\bm{x}) &:=A_{\varphi}(\bm{x})|_{\gamma}, \\
q(\bm{x}) &:=\partial_{\bm{n}}A_{\varphi}(\bm{x})|_{\gamma}
\end{split}
\label{eq:boundary_unknowns}
\end{equation}
represent the azimuthal component and its normal derivative on the boundary, respectively. The normal derivative is evaluated as
\begin{equation}
\frac{\partial}{\partial n} \equiv n_r \frac{\partial}{\partial r} + n_z \frac{\partial}{\partial z}
\label{eq:normal_derivative_rz}
\end{equation}
In \eqref{eq:axisym_bie_air} and \eqref{eq:axisym_bie_conductor}, $\mathcal{G}$ and $\mathcal{G}_k$ denote the axisymmetric Laplace and Helmholtz kernels, respectively, whose explicit expressions are provided in Appendix~C.

\subsection{Galerkin BIE and discretization}

Multiplying both sides of \eqref{eq:axisym_bie_air} by the test function $\phi_i(\bm{x})$ and integrating with respect to $\bm{x}$ over the boundary $\gamma$ yields the Galerkin BIE:
\begin{equation}
\begin{split}
& \int_{\gamma} \phi_i(\bm{x}) \Big\{ c(\bm{x}) u^{(1)}(\bm{x}) \\
& \quad + \int_{\gamma} \left[\mathcal{G}(\bm{x},\bm{y})q^{(1)}(\bm{y})
- u^{(1)}(\bm{y})\partial_{\bm{n}}\mathcal{G}(\bm{x},\bm{y})\right] r_{\bm{y}} ds_{\bm{y}} \\
& \quad - u^{(e)}(\bm{x}) \Big\} r_{\bm{x}} ds_{\bm{x}} = 0 .
\end{split}
\label{eq:galerkin_bie_air}
\end{equation}

Similarly, applying the Galerkin procedure to \eqref{eq:axisym_bie_conductor} gives
\begin{equation}
\begin{split}
& \int_{\gamma} \phi_i(\bm{x}) \Big\{ [1-c(\bm{x})] u^{(2)}(\bm{x}) \\
& \quad + \int_{\gamma} \left[u^{(2)}(\bm{y})\partial_{\bm{n}}\mathcal{G}_k(\bm{x},\bm{y})-\mathcal{G}_k(\bm{x},\bm{y})q^{(2)}(\bm{y})
\right] r_{\bm{y}} ds_{\bm{y}} \Big\} r_{\bm{x}} ds_{\bm{x}} = 0 .
\end{split}
\label{eq:galerkin_bie_conductor}
\end{equation}

In the Galerkin BEM framework, the free-term coefficient $c(\bm{x})$ is set uniformly to $1/2$ along $\gamma$, since variations at corner points, e.g. $1/4$ for right angles, lie on a set of measure zero and thus do not affect the Lebesgue integral \cite{Folland1999}. This convention simplifies the implementation without compromising accuracy. For the discretization, the unknown boundary fields are approximated as linear combinations of boundary basis functions:
\begin{equation}
\begin{aligned}
u_h^{(l)}(\bm{y}) &= \sum_{j=1}^{N_0} u_{j}^{(l)} \phi_j(\bm{y}), \\
q_h^{(l)}(\bm{y}) &= \sum_{j=1}^{N_0} q_{j}^{(l)} \phi_j(\bm{y}),
\end{aligned}
\qquad l \in \{1, 2\}.
\label{eq:field_expansions}
\end{equation}
Substituting \eqref{eq:field_expansions} into \eqref{eq:galerkin_bie_air} directly leads to a system of algebraic equations involving the discrete boundary operators:
\begin{equation}
\left(\frac{1}{2}\mathbf{M}-\mathbf{K}\right)\bm{u}^{(1)}+\mathbf{V}\bm{q}^{(1)}=\bm{f},
\label{eq:discrete_air_system}
\end{equation}
where
\begin{subequations}
\label{eq:laplace_operator_defs}
\begin{align}
M_{ij} &= \int_{\gamma}\phi_i(\bm{x})\phi_j(\bm{x})r_{\bm{x}}ds_{\bm{x}},
\label{eq:operator_M}
\\
V_{ij} &= \int_{\gamma}\int_{\gamma}\phi_i(\bm{x})\mathcal{G}(\bm{x},\bm{y})\phi_j(\bm{y})r_{\bm{y}}ds_{\bm{y}}r_{\bm{x}}ds_{\bm{x}},
\label{eq:operator_V_laplace}
\\
K_{ij} &= \int_{\gamma}\int_{\gamma}\phi_i(\bm{x})\partial_{\bm{n}}\mathcal{G}(\bm{x},\bm{y})\phi_j(\bm{y})r_{\bm{y}}ds_{\bm{y}}r_{\bm{x}}ds_{\bm{x}},
\label{eq:operator_K_laplace}
\\
f_i &= \int_{\gamma}\phi_i(\bm{x})u^{(e)}(\bm{x})r_{\bm{x}}ds_{\bm{x}} .
\label{eq:rhs_excitation}
\end{align}
\end{subequations}

Similarly, the discretization of \eqref{eq:galerkin_bie_conductor} gives
\begin{equation}
\left(\frac{1}{2}\mathbf{M}+\mathbf{K}_k\right)\bm{u}^{(2)}-\mathbf{V}_k\bm{q}^{(2)}=\mathbf{0},
\label{eq:discrete_conductor_system}
\end{equation}
where $\mathbf{V}_k$ and $\mathbf{K}_k$ are obtained by replacing $\mathcal{G}$ with $\mathcal{G}_k$ in \eqref{eq:operator_V_laplace} and \eqref{eq:operator_K_laplace}, respectively.

Additionally, the following interface conditions hold on $\gamma$:
\begin{subequations}
\label{eq:interface_conditions}
\begin{align}
u^{(1)} &= u^{(2)},
\label{eq:interface_u}
\\
q^{(2)} &= \mu_r q^{(1)}+(\mu_{r}-1)n_r \frac{u^{(1)}}{r},
\label{eq:interface_q}
\end{align}
\end{subequations}
where $\mu_r=\mu/\mu_0$ is the relative permeability of the conductive region. Consequently, from \eqref{eq:discrete_air_system}, \eqref{eq:discrete_conductor_system}, \eqref{eq:interface_u}, and \eqref{eq:interface_q}, we obtain
\begin{equation}
\begin{bmatrix}
\frac{1}{2}\mathbf{M}-\mathbf{K} & \mathbf{V} \\
\frac{1}{2}\mathbf{M}+\mathbf{K}_k-(\mu_r-1)\mathbf{V}_s & -\mu_r\mathbf{V}_k
\end{bmatrix}
\begin{bmatrix}
\bm{u}^{(1)} \\
\bm{q}^{(1)}
\end{bmatrix}
=
\begin{bmatrix}
\bm{f} \\
\bm{0}
\end{bmatrix},
\label{eq:final_linear_system}
\end{equation}
where
\begin{equation}
V_{s,ij}=\int_{\gamma}\int_{\gamma}\phi_i(\bm{x})n_r(\bm{y})\mathcal{G}_k(\bm{x},\bm{y})\phi_j(\bm{y})ds_{\bm{y}}r_{\bm{x}}ds_{\bm{x}} .
\label{eq:Vs_operator}
\end{equation}

Subsequently, the impedance variation $\Delta Z$ of the induction coil due to the presence of the conductive object is evaluated via Auld's formula \cite{Auld_Moulder1999}:
\begin{equation}
\Delta Z = -\frac{2\pi i\omega}{\mu_0 I^2} \int_{\gamma} \left( u^{(e)} q^{(1)} - u^{(1)} q^{(e)} \right) r \, ds ,
\label{eq:auld_impedance}
\end{equation}
where $u^{(e)}$ and $q^{(e)}$ denote the excitation vector potential and its normal derivative, while $u^{(1)}$ and $q^{(1)}$ represent the total fields in the non-conductive region obtained from the BEM solution. 

\rev{Since the integrands in \eqref{eq:auld_impedance} are regular, the impedance variation is straightforwardly computed using GL quadrature over the boundary mesh.}

\section{Coordinate Transformation for 2D Singular Integrals}

In the Galerkin BEM framework, evaluating element integrals of the form
\rev{
\begin{equation}
I = \int_0^1 \int_0^1 K(\xi,\eta)\,S(\xi,\eta)\,\mathrm{d}\xi\,\mathrm{d}\eta
\label{eq:galerkin_element_integral}
\end{equation}
requires special treatment when the kernel $K$ is singular. Here, $S(\xi,\eta)$ absorbs the regular shape functions and parametric Jacobians.} The 2D boundary kernels typically exhibit either logarithmic singularities of the form $\ln|x(\xi)-y(\eta)|$ or Cauchy singularities $(x(\xi)-y(\eta))^{-1}$ as source and field points approach coincidence \cite{Guiggiani1991}, rendering direct tensor-product GL quadrature unreliable.

A natural starting point is the Sauter--Schwab quadrature (SSQ)~\cite{Sauter_Schwab2011,Sauter_Schwab1997}, which provides a systematic framework for regularizing singular boundary integrals through coordinate transformations. Originally developed for 4D surface--surface integrals in 3D BEM, the SSQ classifies element pairs into three configurations: identical panels, common edge, and common vertex. Each configuration is addressed by a tailored sequence of Duffy-type mappings \cite{Duffy1982}.

Although the present problem is intrinsically 2D, involving curve--curve interactions, the core principle of the SSQ extends naturally: singular behavior is isolated by subdividing the parameter domain and applying suitable coordinate transformations. In 2D, only two singular configurations arise:
\emph{Coincident elements}, whose parameter domains overlap; and
\emph{Touching elements}, which share a common endpoint. The coordinate transformations employed in this work are summarized below.

\paragraph{Coincident elements}
The parameter square $[0,1]^2$ is divided along the diagonal into two triangles. This partition isolates the singular line $\eta=\xi$, after which a Duffy-type mapping is applied to each region. For the upper triangle $(\eta>\xi)$, we use
\begin{equation}
\xi_1 = u, \qquad \eta_1 = u + (1-u)v, \qquad (u,v)\in[0,1]^2,
\end{equation}
whereas for the lower triangle $(\xi>\eta)$ we use
\begin{equation}
\xi_2 = u + (1-u)v, \qquad \eta_2 = u.
\end{equation}
In both mappings, the Jacobian contributes a factor $|\det J|=1-u$, which transfers the diagonal singularity to the coordinate edge $v=0$.

{\revcolor
For Cauchy-singular kernels, the two mapped triangular contributions must not be evaluated as independent numerical integrals. After both triangles are mapped to the same $(u,v)$-square, their sum is used as the quadrature integrand:
\begin{equation}
I = \int_0^1 \int_0^1 \Big[ K(\xi_1, \eta_1) S(\xi_1, \eta_1) + K(\xi_2, \eta_2) S(\xi_2, \eta_2) \Big] (1-u) \,\mathrm{d}u\,\mathrm{d}v,
\end{equation}
where $(\xi_1,\eta_1)=(u,u+(1-u)v)$ and
$(\xi_2,\eta_2)=(u+(1-u)v,u)$ are the two mirror images of the same point $(u,v)$ in the upper and lower triangles.
For logarithmic kernels, the Jacobian factor $1-u$ is sufficient to regularize the diagonal singularity. For Cauchy-singular kernels, by contrast, regularization requires the antisymmetric recombination of the two mirror triangular contributions before quadrature is applied. The resulting integrand is bounded on the unit square, as proved in Appendix~A. A graded substitution $v=w^p$ with $p>1$ may also be applied after this recombination to concentrate quadrature nodes near $v=0$.
}

\paragraph{Touching elements}
When two elements share a common endpoint, the kernel singularity occurs at a corner of the parameter square $[0,1]^2$. We first consider the canonical touching orientation $x(1)=y(0)$, for which the singular corner is $(\xi,\eta)=(1,0)$. {\revcolor The square is split along the line $\xi+\eta=1$ into the two triangles}
\[
{\revcolor
T_1=\{(\xi,\eta)\in[0,1]^2:\xi+\eta\le 1\},\qquad
T_2=\{(\xi,\eta)\in[0,1]^2:\xi+\eta\ge 1\}.
}
\]
\rev{The two triangles are mapped to the unit square by}
\rev{
\begin{equation}
\xi_1 = 1-u, \qquad \eta_1 = uv, \qquad (u,v)\in[0,1]^2,
\label{eq:touching_map_1}
\end{equation}}
and
\rev{
\begin{equation}
\xi_2 = 1-uv, \qquad \eta_2 = u, \qquad (u,v)\in[0,1]^2.
\label{eq:touching_map_2}
\end{equation}}
Both mappings have Jacobian
\[
{\revcolor
\mathrm{d}\xi\,\mathrm{d}\eta = u\,\mathrm{d}u\,\mathrm{d}v.
}
\]
Accordingly, the transformed contributions are
\rev{
\begin{equation}
I_1 = \int_0^1 \int_0^1 K(1-u, uv)\,S(1-u, uv)\,u\,\mathrm{d}u\,\mathrm{d}v,
\end{equation}}
and
\rev{
\begin{equation}
I_2 = \int_0^1 \int_0^1 K(1-uv, u)\,S(1-uv, u)\,u\,\mathrm{d}u\,\mathrm{d}v.
\end{equation}}

For logarithmic kernels, the regularization follows from the local behavior of the distance to the touching point. Applying \eqref{eq:touching_map_1} gives
\[
{\revcolor
|x(1-u)-y(uv)|
\approx
u\,\bigl|x'(1)+v\,y'(0)\bigr|,
}
\]
whereas applying \eqref{eq:touching_map_2} gives
\[
{\revcolor
|x(1-uv)-y(u)|
\approx
u\,\bigl|v\,x'(1)+y'(0)\bigr|.
}
\]
Thus the logarithmic singularity is converted into a $\ln u$ factor, which is regularized by the Jacobian factor $u$ in both transformed integrals.

For Cauchy-singular kernels, near $(\xi,\eta)=(1,0)$ the integrand can be written as
\rev{
\begin{equation}
K(\xi,\eta)\,S(\xi,\eta) = \frac{A(\xi,\eta)}{(1-\xi)+\eta} ,
\end{equation}}
where $A(\xi,\eta)$ is bounded and smooth on $[0,1]^2$. {\revcolor After applying the two mappings, the denominator becomes}
\[
{\revcolor
(1-\xi_i)+\eta_i = u(1+v),\qquad i=1,2.
}
\]
{\revcolor Therefore
\[
\frac{A(\xi_i,\eta_i)}{(1-\xi_i)+\eta_i}\,u
=
\frac{A(\xi_i,\eta_i)}{1+v},
\qquad i=1,2,
\]
which is bounded over the unit square. Thus, the singularities in both $I_1$ and $I_2$ are completely eliminated.}
The full touching-element contribution is then
\[
{\revcolor
I=I_1+I_2.
}
\]
\rev{For the canonical orientation $x(1)=y(0)$, this approach removes the endpoint singularity at $(1,0)$ and yields regular double integrals over the unit square. The opposite orientation $x(0)=y(1)$, whose singular corner is $(\xi,\eta)=(0,1)$, is reduced to the same case by reversing both local parameters,
\[
    \xi \to 1-\xi,\qquad \eta \to 1-\eta .
\]
Hence the same regularized rules for $I_1$ and $I_2$ are reused without any additional transformation formula.}

\section{Numerical Validation of Axisymmetric Galerkin BEM}
To evaluate the accuracy and convergence characteristics of the proposed formulation, three representative axisymmetric configurations are analyzed: a cylindrical tube, a conical tube, and a spherical shell. In each case, the impedance change of an axisymmetric eddy-current coil positioned in proximity to a conducting structure is computed. A common coil geometry is employed across all simulations (Table~\ref{tab:coil_common}) \rev{so that the validation of each configuration focuses solely on the boundary discretization and the proposed integration algorithms, rather than on variations in the source field.} The coil's source potential $A_\varphi^{(e)}$ and its normal derivative are obtained through analytical expressions provided in Appendix B. Impedance variations are normalized and reported as $\Delta R/X_0$ and $\Delta X/X_0$, where $X_0=\omega L_0$ denotes the reactance of the isolated coil with self-inductance $L_0$ (Table~\ref{tab:coil_common}).

\begin{table}[H]
\centering
\caption{Common coil parameters used in all validation cases.}
\label{tab:coil_common}
\small
\setlength{\tabcolsep}{8pt}
\renewcommand{\arraystretch}{1.15}
\begin{tabular}{@{} l c l @{}}
\toprule\addlinespace[.25ex]\toprule
\textbf{Parameter} & \textbf{Symbol} & \textbf{Value} \\
\midrule
Inner radius (m)           & $r_1$   & \num{0.007} \\
Outer radius (m)           & $r_2$   & \num{0.0085} \\
Half-height (m)            & $h$     & \num{0.002} \\
Turns                      & $N$     & \num{500} \\
Baseline inductance (mH)   & $L_0$   & \num{4.7405622} \\
\bottomrule\addlinespace[.25ex]\bottomrule
\end{tabular}
\end{table}

{\revcolor
The FEM reference solutions for all three validation cases were computed with COMSOL Multiphysics 6.1. Infinite element domains were applied to the exterior region to avoid artificial truncation of the open boundary. The conducting regions were meshed according to a skin-depth criterion at the highest frequency of \(100\,\mathrm{kHz}\), with approximately ten elements per skin depth in each case. The case-dependent FEM mesh settings, system sizes, and solution times are summarized in Table~\ref{tab:fem_reference_settings}. The proposed Galerkin BEM was implemented separately in Julia. Both the COMSOL simulations and the Julia BEM computations were performed on the same laptop equipped with an AMD Ryzen AI 9 HX 370 processor (2.0\,GHz) and 32\,GB of RAM.
}
\begin{table}[H]
\centering
\caption{\rev{FEM reference settings used in the validation cases. The mesh size denotes the maximum element size in the conducting region.}}
\label{tab:fem_reference_settings}
\small
\setlength{\tabcolsep}{6pt}
\renewcommand{\arraystretch}{1.15}
{\revcolor
\begin{tabular}{@{} l c c c c @{}}
\toprule\addlinespace[.25ex]\toprule
\textbf{Case} &
\textbf{Material} &
\textbf{Mesh size} &
\textbf{DoFs} &
\textbf{Solution time} \\
\midrule
Cylindrical tube &
S30400 &
$130\,\mu\mathrm{m}$ &
\num{74665} &
$4\,\mathrm{s}$ \\
Conical tube &
7075--T6 &
$40\,\mu\mathrm{m}$ &
\num{263768} &
$6\,\mathrm{s}$ \\
Spherical shell &
C96400 &
$90\,\mu\mathrm{m}$ &
\num{106163} &
$5\,\mathrm{s}$ \\
\bottomrule\addlinespace[.25ex]\bottomrule
\end{tabular}
}
\end{table}

\subsection{A coil placed inside a S30400 stainless steel tube}
We consider a coaxial eddy-current inspection configuration in which a cylindrical coil is placed inside a stainless-steel tube (AISI 304, UNS S30400). The tube has inner and outer radii $a_1$, $a_2$, and axial length $l$, while the coil has inner and outer radii $r_1$ and $r_2$ and half-height $h$.

\begin{figure}[t]
\centering
\begin{tikzpicture}[
    scale=6.0, >=Latex, line cap=round, line join=round,
    every node/.style={font=\small,align=center},
    axis/.style={black, very thick},
    body/.style={black, very thick},
    coil/.style={black, very thick, fill=black!10},
    dim/.style={black, thin},
    guide/.style={densely dashed, thin}
]

\def\aone{0.82}    
\def\atwo{0.96}    
\def\Lh  {0.60}    
\def\rone{0.50}    
\def\rtwo{0.66}    
\def\hh  {0.18}    
\def\zO  {0.25}    
\def\zA {-0.38}    
\def\xTwoH{0.70}   

\draw[->,axis] (-0.06,0) -- (1.1,0) node[below] {$r$};
\draw[->,axis] (0,-0.75) -- (0, 0.75) node[left] {$z$};
\node[below left] at (0,0) {$0$};

\draw[body] (\aone,\Lh) -- (\atwo,\Lh) -- (\atwo,-\Lh) -- (\aone,-\Lh) -- cycle;
\node at ({(\aone+\atwo)/2},0.09) {$\sigma$\\$\mu_r$};

\draw[guide] (0,\Lh) -- (\aone,\Lh);
\draw[guide] (0,-\Lh) -- (\aone,-\Lh);
\node[left]  at (0,\Lh)  {$+\,l/2$};
\node[left]  at (0,-\Lh) {$-\,l/2$};

\draw[coil] (\rone,\zO-\hh) rectangle (\rtwo,\zO+\hh);

\draw[guide] (0,\zO) -- (\rone,\zO);
\node[left] at (0,\zO) {$z_0$};

\pgfmathsetmacro{\rc}{(\rone+\rtwo)/2}
\def\rad{0.04}
\pgfmathsetmacro{\s}{\rad/sqrt(2)}
\draw (\rc,\zO) circle (\rad);
\draw (\rc-\s,\zO-\s) -- (\rc+\s,\zO+\s);
\draw (\rc-\s,\zO+\s) -- (\rc+\s,\zO-\s);

\draw[guide] (\rtwo,\zO+\hh) -- (\xTwoH,\zO+\hh);
\draw[guide] (\rtwo,\zO-\hh) -- (\xTwoH,\zO-\hh);
\draw[<->,dim] (\xTwoH,\zO-\hh) -- node[right] {$2h$} (\xTwoH,\zO+\hh);

\draw[guide] (\rone,\zO+\hh) -- (\rone,0);
\draw[guide] (\rtwo,\zO+\hh) -- (\rtwo,0);
\node[below] at (\rone,0) {$r_1$};
\node[below] at (\rtwo,0) {$r_2$};

\draw[guide] (\aone,0) -- (\aone,\zA);
\draw[guide] (\atwo,0) -- (\atwo,\zA-0.10);
\draw[<->,dim] (0,\zA) -- node[below] {$a_1$} (\aone,\zA);
\draw[<->,dim] (0,\zA-0.10) -- node[below] {$a_2$} (\atwo,\zA-0.10);

\end{tikzpicture}
\caption{\rev{Meridional cross-section of a conductive tube and a coaxial excitation coil.}}
\label{fig:geom-304}
\end{figure}

\begin{figure}[H]
\centering
\includegraphics[width=0.95\linewidth]{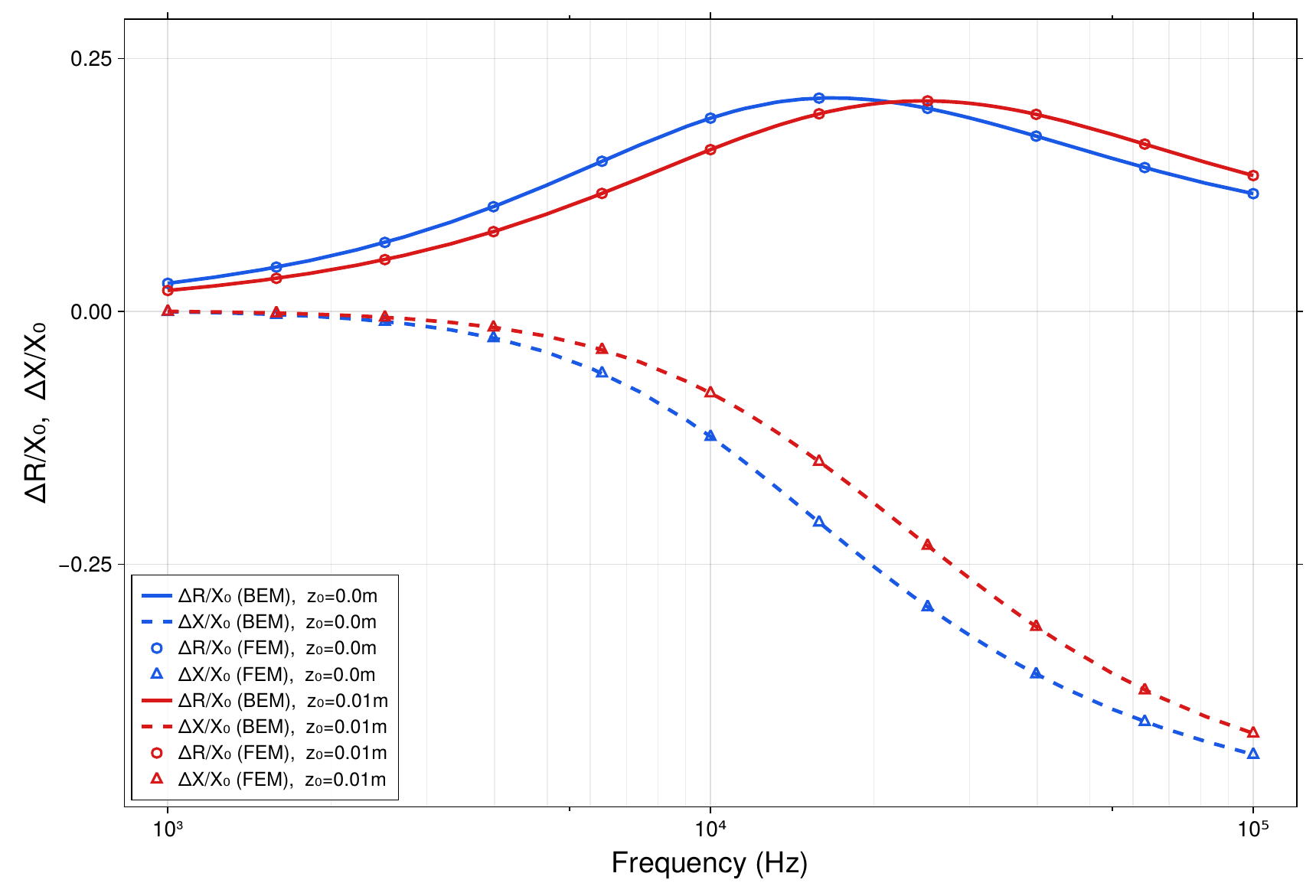}
\caption{Comparison of $\Delta R/X_0$ and $\Delta X/X_0$ obtained by the Galerkin BEM and FEM for the S30400 tube at two lift-off positions. \rev{The BEM predictions are shown as continuous curves, while the FEM reference solutions are denoted by discrete markers.}}
\label{fig:CaseA_sweep}
\end{figure}

The conductive tube is characterized by relative permeability $\mu_r$ and conductivity $\sigma$, giving rise to the complex wavenumber $k = \sqrt{-i\omega\mu_0 \mu_r \sigma}$. The tube surface is discretized with linear boundary elements, and all singular integrals are evaluated using the algorithm in \textsection 3. \rev{Figure~\ref{fig:CaseA_sweep} plots the normalized impedance spectra ($\Delta R/X_0$ and $\Delta X/X_0$) over a given frequency range for two axial coil positions. The BEM predictions show excellent agreement with the FEM references (settings in Table~\ref{tab:fem_reference_settings})}.

\begin{table}[H]
\centering
\caption{Geometric and material parameters for the S30400 stainless-steel tube.}
\label{tab:caseA_tube}
\small
\setlength{\tabcolsep}{8pt}
\renewcommand{\arraystretch}{1.15}
\begin{tabular}{@{} l c l @{}}
\toprule\addlinespace[.25ex]\toprule
\textbf{Parameter} & \textbf{Symbol} & \textbf{Value} \\
\midrule
Inner tube radius (m)      & $a_1$   & \num{0.009} \\
Outer tube radius (m)      & $a_2$   & \num{0.011} \\
Tube length (m)            & $l$     & \num{0.024} \\
Relative permeability      & $\mu_r$ & \num{1.021} \\
Conductivity (MS/m)        & $\sigma$& \num{1.37} \\
\bottomrule\addlinespace[.25ex]\bottomrule
\end{tabular}
\end{table}

\subsection{A coil located inside a 7075-T6 aluminum conical tube}
We next consider a conical tube fabricated from 7075-T6 aluminum, with inner and outer radii that vary linearly along the $z$-axis. The smaller end of the tube is located at $z=-l/2$ and the larger end at $z=+l/2$, maintaining a uniform wall thickness throughout. The geometry is illustrated in Fig.~\ref{fig:geom-conical}. The coil parameters remain identical to those of the preceding cylindrical case.

\begin{figure}[H]
\centering
\begin{tikzpicture}[
    scale=6.0, >=Latex, line cap=round, line join=round,
    every node/.style={font=\small,align=center},
    axis/.style={black, very thick},
    body/.style={black, very thick},
    coil/.style={black, very thick, fill=black!10},
    dim/.style={black, thin},
    guide/.style={densely dashed, thin}
]

\def\Lh  {0.60}     
\def\aoneBot{0.48}  
\def\aoneTop{0.62}  
\def\thick  {0.10}  
\def\atwoBot{\aoneBot+\thick} 
\def\atwoTop{\aoneTop+\thick} 

\def\rone{0.34}     
\def\rtwo{0.44}     
\def\hh  {0.14}     
\def\zO  {0.18}     

\def\xTwoH {0.47}

\def\deltaTop{0.10}   
\def\deltaBot{0.10}   
\pgfmathsetmacro{\yTopBelow}{\Lh - \deltaTop}
\pgfmathsetmacro{\yBotAbove}{-\Lh + \deltaBot}

\draw[->,axis] (-0.06,0) -- (0.80,0) node[below] {$r$};
\draw[->,axis] (0,-0.70) -- (0,0.70) node[left] {$z$};
\node[below left] at (0,0) {$0$};

\draw[guide] (0,\Lh) -- (\aoneTop,\Lh);
\draw[guide] (0,-\Lh) -- (\aoneBot,-\Lh);
\node[left]  at (0,\Lh)  {$+\,l/2$};
\node[left]  at (0,-\Lh) {$-\,l/2$};

\draw[body] (\aoneBot,-\Lh) -- (\aoneTop,\Lh);   
\draw[body] (\atwoBot,-\Lh) -- (\atwoTop,\Lh);   
\draw[body] (\aoneTop,\Lh) -- (\atwoTop,\Lh);    
\draw[body] (\aoneBot,-\Lh) -- (\atwoBot,-\Lh);  

\draw[coil] (\rone,\zO-\hh) rectangle (\rtwo,\zO+\hh);

\draw[guide] (0,\zO) -- (\rone,\zO);
\node[left] at (0,\zO) {$z_0$};

\pgfmathsetmacro{\rc}{(\rone+\rtwo)/2}
\def\rad{0.035}
\pgfmathsetmacro{\s}{\rad/sqrt(2)}
\draw (\rc,\zO) circle (\rad);
\draw (\rc-\s,\zO-\s) -- (\rc+\s,\zO+\s);
\draw (\rc-\s,\zO+\s) -- (\rc+\s,\zO-\s);

\draw[guide] (\rtwo,\zO+\hh) -- (\xTwoH,\zO+\hh);
\draw[guide] (\rtwo,\zO-\hh) -- (\xTwoH,\zO-\hh);
\draw[<->,dim] (\xTwoH,\zO-\hh) -- node[right] {$2h$} (\xTwoH,\zO+\hh);

\draw[guide] (\rone,\zO+\hh) -- (\rone,0);
\draw[guide] (\rtwo,\zO+\hh) -- (\rtwo,0);
\node[below] at (\rone,0) {$r_1$};
\node[below] at (\rtwo,0) {$r_2$};

\draw[guide] (\aoneTop,\Lh) -- (\aoneTop,\yTopBelow);
\draw[guide] (\atwoTop,\Lh) -- (\atwoTop,\yTopBelow-0.10);
\draw[<->,dim] (0,\yTopBelow) -- node[above] {$a_3$} (\aoneTop,\yTopBelow);
\draw[<->,dim] (0,\yTopBelow-0.10) -- node[above] {$a_4$} (\atwoTop,\yTopBelow-0.10);

\draw[guide] (\aoneBot,-\Lh) -- (\aoneBot,\yBotAbove);
\draw[guide] (\atwoBot,-\Lh) -- (\atwoBot,\yBotAbove+0.10);
\draw[<->,dim] (0,\yBotAbove) -- node[below] {$a_1$} (\aoneBot,\yBotAbove);
\draw[<->,dim] (0,\yBotAbove+0.10) -- node[below] {$a_2$} (\atwoBot,\yBotAbove+0.10);

\pgfmathsetmacro{\rinLab }{\aoneBot + (\aoneTop-\aoneBot)*((\zO+\Lh)/(2*\Lh))}
\pgfmathsetmacro{\routLab}{\atwoBot + (\atwoTop-\atwoBot)*((\zO+\Lh)/(2*\Lh))}
\def\f{0.41}
\def\eps{0.005}
\pgfmathsetmacro{\rpos}{\rinLab + \f*(\routLab-\rinLab) - \eps}
\node[anchor=east, inner sep=1pt] at (\rpos, \zO+0.035) {$\sigma$};
\node[anchor=east, inner sep=1pt] at (\rpos, \zO-0.035) {$\mu_r$};

\end{tikzpicture}

\caption{\rev{Meridional cross-section of a conical conductive tube and a coaxial excitation coil.}}
\label{fig:geom-conical}
\end{figure}

\begin{figure}[H]
\centering
\includegraphics[width=0.95\linewidth]{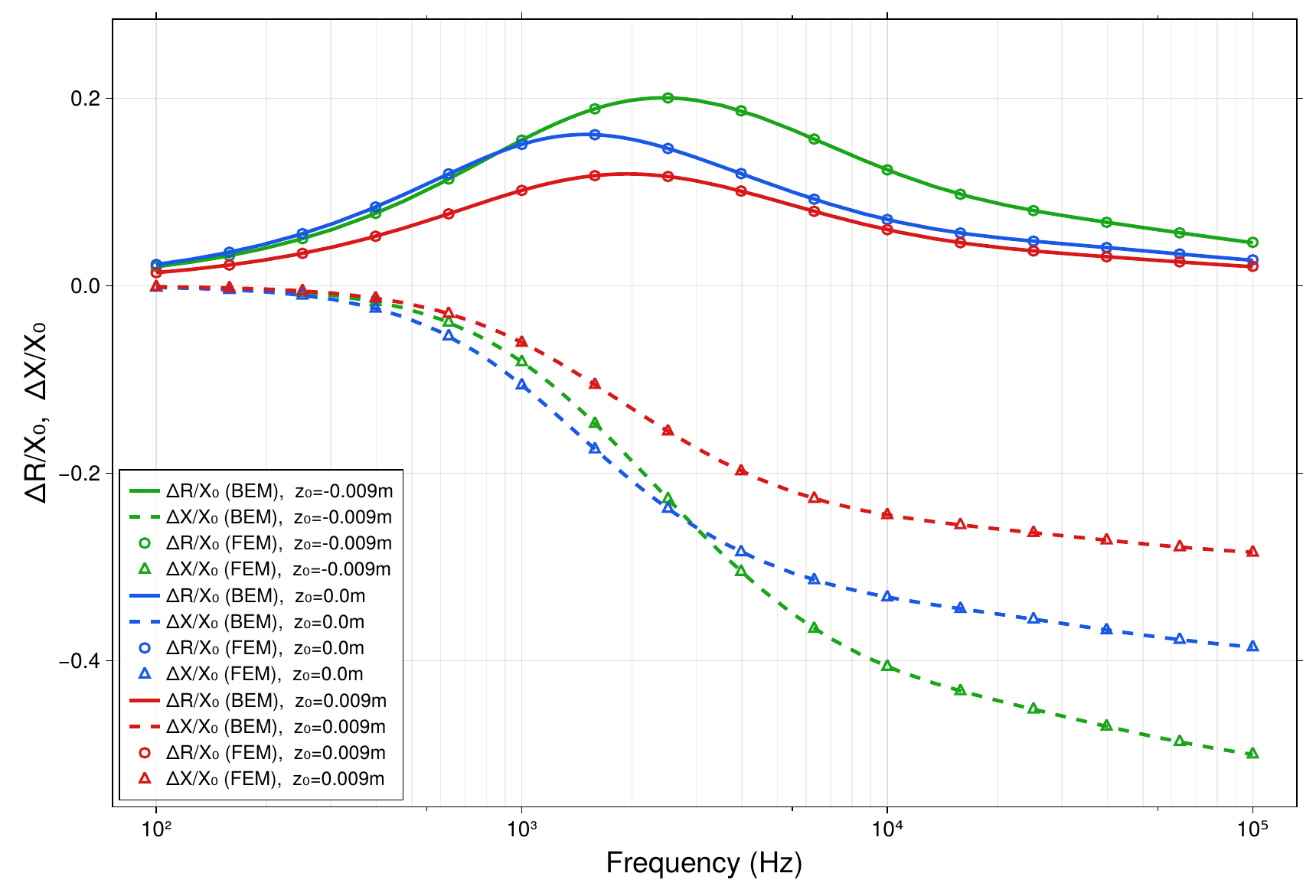}
\caption{Comparison of $\Delta R/X_0$ and $\Delta X/X_0$ obtained by the Galerkin BEM and the FEM reference for the conical 7075--T6 tube at different axial positions of the coil. \rev{The BEM predictions are shown as continuous curves, while the FEM reference solutions are denoted by discrete markers.}}
\label{fig:CaseB_sweep}
\end{figure}

For the conical tube validation, the slanted boundary is discretized with linear elements, applying the singular integration algorithm from \textsection 3. Impedance variations are computed for three axial coil positions over the same frequency range. Figure~\ref{fig:CaseB_sweep} compares the normalized resistance and reactance changes ($\Delta R/X_0$ and $\Delta X/X_0$), obtained via Auld's formula, against the FEM reference data \rev{(settings in Table~\ref{tab:fem_reference_settings})}. Excellent agreement is achieved across all frequencies and coil positions.

\begin{table}[H]
\centering
\caption{Geometric and material parameters for the conical 7075--T6 aluminum tube.}
\label{tab:caseB_cone}
\small
\setlength{\tabcolsep}{8pt}
\renewcommand{\arraystretch}{1.15}
\begin{tabular}{@{} l c l @{}}
\toprule\addlinespace[.25ex]\toprule
\textbf{Parameter} & \textbf{Symbol} & \textbf{Value} \\
\midrule
Bottom inner radius (m)     & $a_1$ & \num{0.0090} \\
Bottom outer radius (m)     & $a_2$ & \num{0.0106} \\
Top inner radius (m)        & $a_3$ & \num{0.0108} \\
Top outer radius (m)        & $a_4$ & \num{0.0124} \\
Tube length (m)             & $l$   & \num{0.022} \\
Relative permeability       & $\mu_r$ & \num{1.0} \\
Conductivity (MS/m)         & $\sigma$& \num{16.93} \\
\bottomrule\addlinespace[.25ex]\bottomrule
\end{tabular}
\end{table}

\subsection{A coil placed near a C96400 copper alloy spherical shell}
The third validation examines a spherical shell fabricated from C96400 copper alloy \cite{Luo2023}, with the coil positioned either external or internal to the shell. This geometry provides a more stringent test of the numerical framework, as the curved boundary demands accurate integration over curved segments and introduces global coupling among disjoint regions. Coil specifications are given in Table~\ref{tab:coil_common}, while shell dimensions and material properties appear in Table~\ref{tab:caseC_sphere}. The configuration is illustrated schematically in Figure~\ref{fig:geom-sphere}.

\begin{figure}[H]
\centering
\includegraphics[width=0.95\linewidth]{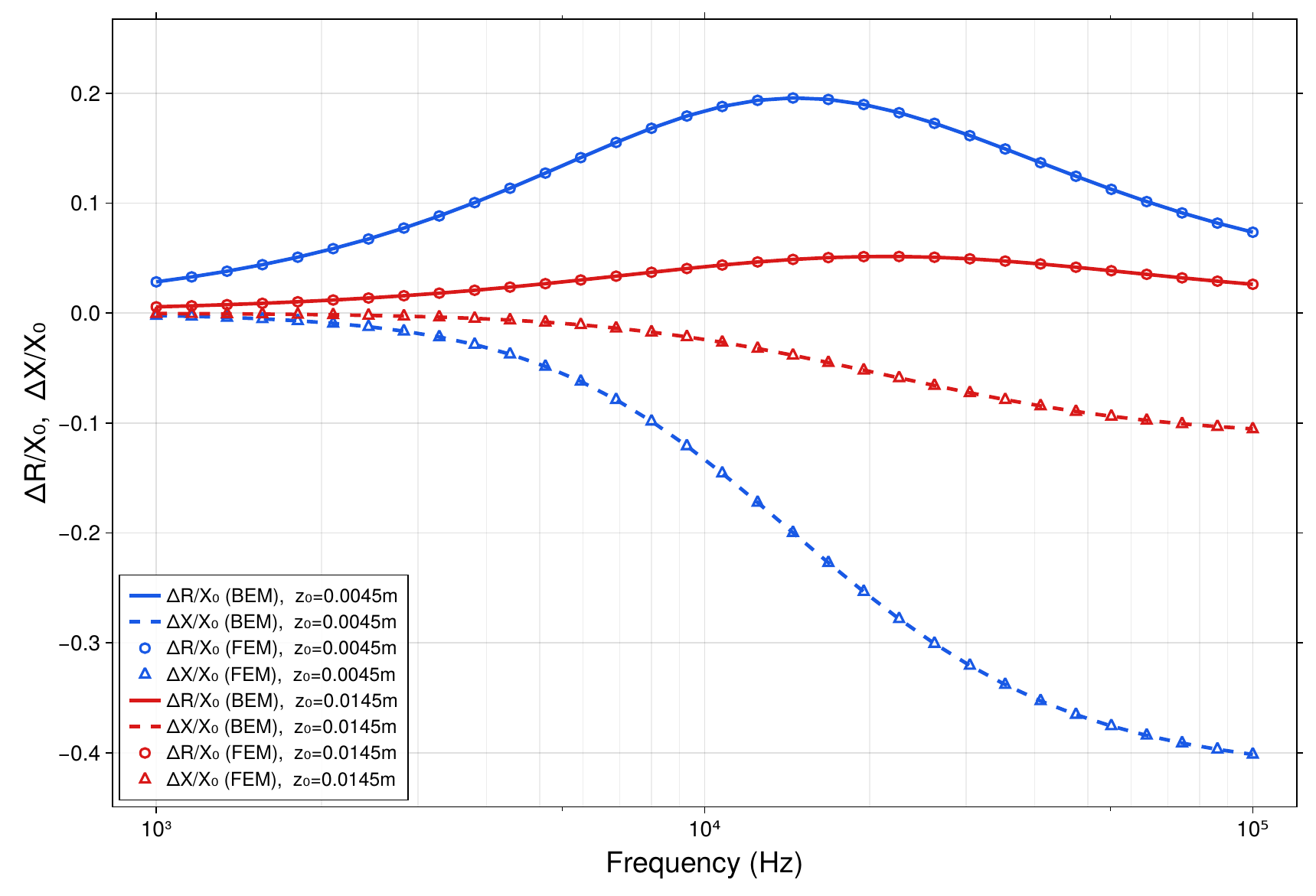}
\caption{Comparison of $\Delta R/X_0$ and $\Delta X/X_0$ obtained by the Galerkin BEM and the FEM reference for the C96400 spherical shell at two axial positions of the coil. \rev{The BEM predictions are shown as continuous curves, while the FEM reference solutions are denoted by discrete markers.}}
\label{fig:CaseC_sweep}
\end{figure}

All singular kernels are evaluated using the algorithm detailed in \textsection 3, ensuring numerical accuracy for every element interaction. Impedance variations are computed at two axial coil positions, and the normalized components $\Delta R/X_0$ and $\Delta X/X_0$, obtained via Auld's formula, are compared with \rev{FEM reference solutions (settings detailed in Table~\ref{tab:fem_reference_settings}).} Figure~\ref{fig:CaseC_sweep} presents the frequency-dependent response; close agreement between boundary-element predictions and reference solutions is observed throughout the parameter range investigated.

\begin{figure}[H]
\centering
\begin{tikzpicture}[
    scale=6.0, >=Latex, line cap=round, line join=round,
    every node/.style={font=\small,align=center},
    axis/.style={black, very thick},
    body/.style={black, very thick},
    coil/.style={black, very thick, fill=black!10},
    dim/.style={black, thin},
    guide/.style={densely dashed, thin}
]

\def\aone{0.50}        
\def\thk {0.08}        
\def\atwo{\aone+\thk}  

\def\rone{0.2}
\def\rtwo{0.3}
\def\hh  {0.08}        
\def\zO  {0.15}        
\def\xTwoH{0.35}       

\draw[->,axis] (-0.04,0) -- (0.75,0) node[below] {$r$};
\draw[->,axis] (0,-0.75) -- (0,0.75) node[left] {$z$};
\node[below left] at (0,0) {$0$};

\begin{scope}
  \clip (0,-0.8) rectangle (0.75,0.8); 
  \draw[body] (0,0) circle[radius=\aone];
  \draw[body] (0,0) circle[radius=\atwo];
\end{scope}

\pgfmathsetmacro{\rinLab }{\aone}
\pgfmathsetmacro{\routLab}{\atwo}

\def\fsig{0.15}   
\def\fmu {0.4}   

\def\dysig{0.040}
\def\dymu {0.040}

\pgfmathsetmacro{\rpossig}{\rinLab + \fsig*(\routLab-\rinLab)}
\pgfmathsetmacro{\rposmu }{\rinLab + \fmu *(\routLab-\rinLab)}

\node[anchor=center, inner sep=0pt] at (\rpossig, \zO+\dysig) {$\sigma$};
\node[anchor=center, inner sep=0pt] at (\rposmu , \zO-\dymu ) {$\mu_r$};

\draw[coil] (\rone,\zO-\hh) rectangle (\rtwo,\zO+\hh);
\draw[guide] (0,\zO) -- (\rone,\zO);   
\node[left] at (0,\zO) {$z_0$};

\pgfmathsetmacro{\rc}{(\rone+\rtwo)/2}
\def\rad{0.030}
\pgfmathsetmacro{\s}{\rad/sqrt(2)}
\draw (\rc,\zO) circle (\rad);
\draw (\rc-\s,\zO-\s) -- (\rc+\s,\zO+\s);
\draw (\rc-\s,\zO+\s) -- (\rc+\s,\zO-\s);

\draw[guide] (\rtwo,\zO+\hh) -- (\xTwoH,\zO+\hh);
\draw[guide] (\rtwo,\zO-\hh) -- (\xTwoH,\zO-\hh);
\draw[<->,dim] (\xTwoH,\zO-\hh) -- node[right] {$2h$} (\xTwoH,\zO+\hh);

\draw[guide] (\rone,\zO+\hh) -- (\rone,0);
\draw[guide] (\rtwo,\zO+\hh) -- (\rtwo,0);
\node[below] at (\rone,0) {$r_1$};
\node[below] at (\rtwo,0) {$r_2$};
\node[below] at (\aone-0.04,0) {$a_1$};
\node[below] at (\atwo+0.04,0) {$a_2$};

\end{tikzpicture}
\caption{\rev{Meridional cross-section of a spherical shell and a coaxial excitation coil.}}
\label{fig:geom-sphere}
\end{figure}

\begin{table}[H]
\centering
\caption{Geometric and material parameters for the C96400 spherical shell.}
\label{tab:caseC_sphere}
\small
\setlength{\tabcolsep}{8pt}
\renewcommand{\arraystretch}{1.15}
\begin{tabular}{@{} l c l @{}}
\toprule\addlinespace[.25ex]\toprule
\textbf{Parameter} & \textbf{Symbol} & \textbf{Value} \\
\midrule
Inner radius (m)            & $a_1$   & \num{0.011} \\
Outer radius (m)            & $a_2$   & \num{0.012} \\
Relative permeability       & $\mu_r$ & \num{1.0} \\
Conductivity (MS/m)         & $\sigma$& \num{2.9} \\
\bottomrule\addlinespace[.25ex]\bottomrule
\end{tabular}
\end{table}

Table~\ref{tab:caseC_errors_timing} summarizes the relative errors and average CPU time per frequency point for both the $\mathrm{P}_1$ (linear) and $\mathrm{P}_2$ (quadratic isoparametric) Galerkin BEM formulations applied to the spherical shell. \rev{The mesh parameter $n_s$ denotes the number of boundary elements on each semicircular meridional arc; the total number of boundary elements is therefore $2n_s$.} Relative errors for the $\mathrm{P}_1$ formulation are defined as
\begin{equation}
\begin{split}
& \epsilon^{\mathrm{P}_1}_R = 
\frac{\left\|(\Delta R/X_0)_{\mathrm{P}_1}
-(\Delta R/X_0)_{\mathrm{ref}}\right\|_{2}}{\left\|(\Delta R/X_0)_{\mathrm{ref}}\right\|_{2}}, \\
& \epsilon^{\mathrm{P}_1}_X = 
\frac{\left\|(\Delta X/X_0)_{\mathrm{P}_1}
-(\Delta X/X_0)_{\mathrm{ref}}\right\|_{2}}{\left\|(\Delta X/X_0)_{\mathrm{ref}}\right\|_{2}},
\end{split}
\end{equation}
where $(\Delta R/X_0)_{\mathrm{ref}}$ and $(\Delta X/X_0)_{\mathrm{ref}}$ denote \rev{the BEM reference values computed with the $\mathrm{P}_2$ formulation on a refined mesh ($n_s = 1280$)}. The quantities $\epsilon^{\mathrm{P}_2}_R$ and $\epsilon^{\mathrm{P}_2}_X$ are defined analogously for the $\mathrm{P}_2$ formulation. The entries $t_{\mathrm{P}_1}$ and $t_{\mathrm{P}_2}$ report the average CPU time in seconds required to solve a single frequency point after pre-compilation.

{\revcolor
Both formulations converge monotonically toward the fine-BEM reference under mesh refinement. The singular quadrature used the graded substitution \(v=w^p\) with \(p=3\). More importantly, the $\mathrm{P}_2$ formulation reaches the sub-$0.01\%$ error level already at $n_s=80$, with average CPU times of only $0.14\,\mathrm{s}$ per frequency point. On the same hardware platform, the COMSOL reference model requires about $5\,\mathrm{s}$ per frequency point. Thus, at this practically sufficient accuracy level, the proposed $\mathrm{P}_2$ BEM is about 35 times faster than the FEM benchmark. Such a speed advantage is important in conductivity or permeability inversion, coil-position scans, frequency sweeps, and geometry optimization, where the forward model must be evaluated many times for different trial parameters.
}

\begin{table}[H]
\centering
\caption{\rev{Self-convergence relative errors with respect to a BEM reference ($\mathrm{P}_2$, $n_s = 1280$)} and average CPU time per frequency point (Spherical Shell, $z_0 = 0.0045$ m).}
\label{tab:caseC_errors_timing}
\small
\setlength{\tabcolsep}{8pt}
\renewcommand{\arraystretch}{1.15}
\begin{tabular}{@{} c
                  c c c
                  c c c @{}}
\toprule \addlinespace[.25ex] \toprule
\textbf{$n_s$} &
\textbf{$\epsilon^{\mathrm{P1}}_R$ [\%]} &
\textbf{$\epsilon^{\mathrm{P1}}_X$ [\%]} &
\textbf{$t_{\mathrm{P1}}$ [s]} &
\textbf{$\epsilon^{\mathrm{P2}}_R$ [\%]} &
\textbf{$\epsilon^{\mathrm{P2}}_X$ [\%]} &
\textbf{$t_{\mathrm{P2}}$ [s]} \\
\midrule
40  & \revcolor 0.211   & \revcolor 0.168   & \revcolor 0.03 & \revcolor 0.00859  & \revcolor 0.0119  & \revcolor 0.06 \\
80  & \revcolor 0.0513  & \revcolor 0.0395  & \revcolor 0.07 & \revcolor 0.00361  & \revcolor 0.00502 & \revcolor 0.14 \\
160 & \revcolor 0.0124  & \revcolor 0.00932 & \revcolor 0.22 & \revcolor 0.00130  & \revcolor 0.00181 & \revcolor 0.40 \\
320 & \revcolor 0.00298 & \revcolor 0.00221 & \revcolor 0.80 & \revcolor 0.000415 & \revcolor 0.000576& \revcolor 1.28 \\
\bottomrule \addlinespace[.25ex]
\bottomrule
\end{tabular}
\end{table}

\section{Conclusion}
A general, high-accuracy Galerkin BEM for axisymmetric eddy-current analysis has been developed. \rev{The principal technical contribution is a unified coordinate-transformation framework that regularizes both logarithmic and Cauchy singularities arising from coincident and touching Galerkin element pairs. Consequently, the regularized integrands allow direct GL quadrature, bypassing any kernel- or order-specific analytical extractions.} Benchmark problems involving cylindrical, conical, and spherical conductors validate the method's accuracy and robustness. Linear elements exhibit classical second-order convergence, while quadratic isoparametric elements achieve relative errors on the order of \rev{$10^{-4}$} even on coarse meshes. \rev{These results confirm that the proposed regularization enables accurate, flexible, and efficient eddy-current axisymmetric modeling via high-order Galerkin BEM. The regularization approach is compact, broadly applicable, and readily extensible to other 2D Galerkin BEM beyond eddy-current testing.}

\appendix

\section*{Appendix A. Regularization of Cauchy-singular kernels of coincident elements under the symmetric split Duffy mapping}
\setcounter{equation}{0}
\setcounter{table}{0}
\renewcommand{\theequation}{A\arabic{equation}}
\renewcommand{\thetable}{A\arabic{table}}

For coincident boundary elements that share the same smooth parametrization, the Cauchy-singular kernel in equation \eqref{eq:galerkin_element_integral} of \textsection 3 can be written in the standard form
\rev{
\begin{equation}
K\bigl(\xi,\eta\bigr)S\bigl(\xi,\eta\bigr)
= \frac{F(\xi,\eta)}{\xi-\eta},
\label{eq:A_start}
\end{equation}}
where \(F(\xi,\eta)\) is a smooth function on $[0,1]^2$.
Then, the equation \eqref{eq:galerkin_element_integral} can be written in the principal–value sense as
\begin{equation}
I = \operatorname{p.v.}\int_{0}^{1}\int_{0}^{1}
\frac{F(\xi,\eta)}{\xi-\eta}d\xi d\eta.
\label{eq:A_pv}
\end{equation}

To isolate the singular line \(\eta=\xi\), we split the unit square into the upper and lower triangles
\[
T_1 = \{(\xi,\eta)\in[0,1]^2 : \eta>\xi\},\qquad
T_2 = \{(\xi,\eta)\in[0,1]^2 : \xi>\eta\},
\]
so that
\begin{equation}
I = \int_{T_1}\frac{F(\xi,\eta)}{\xi-\eta}\,d\xi d\eta
+ \int_{T_2}\frac{F(\xi,\eta)}{\xi-\eta}\,d\xi d\eta\ = I_1 + I_2.
\end{equation}
For \(T_1\) we use the Duffy mapping
\begin{equation}
\xi = u,\qquad \eta = u + (1-u)v,\qquad (u,v)\in[0,1]^2,
\end{equation}
which maps the unit square in \((u,v)\) onto the upper triangle \(\eta>\xi\).
The Jacobian determinant is
\[
\left\lvert\frac{\partial(\xi,\eta)}{\partial(u,v)}\right\rvert
= 1-u.
\]
Moreover,
\[
\xi - \eta = -(1-u)v.
\]
Thus
\begin{equation}
\begin{aligned}
I_1
&= \int_0^1\!\!\int_0^1
\frac{F\bigl(u,u+(1-u)v\bigr)}{-(1-u)v}\,(1-u)\,du\,dv \\
&= - \int_0^1\!\!\int_0^1
\frac{F\bigl(u,u+(1-u)v\bigr)}{v}\,du\,dv.
\end{aligned}
\label{eq:A_Iplus}
\end{equation}
For \(T_2\) we use the symmetric mapping
\begin{equation}
\xi = u + (1-u)v,\qquad \eta = u,\qquad (u,v)\in[0,1]^2,
\end{equation}
with the same Jacobian \(\lvert J\rvert=1-u\) and
\[
\xi-\eta = (1-u)v.
\]
Hence
\begin{equation}
\begin{aligned}
I_2
&= \int_0^1\!\!\int_0^1
\frac{F\bigl(u+(1-u)v,u\bigr)}{(1-u)v}\,(1-u)\,du\,dv \\
&= \int_0^1\!\!\int_0^1
\frac{F\bigl(u+(1-u)v,u\bigr)}{v}\,du\,dv.
\end{aligned}
\label{eq:A_Iminus}
\end{equation}
Adding \eqref{eq:A_Iplus} and \eqref{eq:A_Iminus} yields the transformed
representation
\begin{equation}
I = \int_0^1\!\!\int_0^1
\frac{F\bigl(u+(1-u)v,u\bigr) - F\bigl(u,u+(1-u)v\bigr)}{v}\,du\,dv.
\label{eq:A_I_uv}
\end{equation}
The apparent singular factor \(1/v\) is now multiplied by a difference of two evaluations of the same smooth function \(F\). We now show that this difference cancels the Cauchy singularity and yields a bounded integrand on the
whole \([0,1]^2\).

For any fixed \(u\in[0,1]\), since \(F(\xi,\eta)\in C^1([0,1]^2)\), a Taylor expansion of \(F(\xi,\eta)\) about the point \((u,u)\) gives, for small \(v\),
\begin{align}
F\bigl(u+(1-u)v,u\bigr)
&= F(u,u) + (1-u)v\,F_\xi^{'}(u,u) + O(v^2),
\\
F\bigl(u,u+(1-u)v\bigr)
&= F(u,u) + (1-u)v\,F_\eta^{'}(u,u) + O(v^2),
\end{align}
Subtracting these two expansions, we obtain
\begin{equation}
\begin{split}
& F\bigl(u+(1-u)v,u\bigr) - F\bigl(u,u+(1-u)v\bigr) \\
& = (1-u)v\,\bigl[F_\xi^{'}(u,u) - F_\eta^{'}(u,u)\bigr] + O(v^2).
\end{split}
\end{equation}
Hence, as \(v\to 0\), we obtain a regular integrand on the closed square \([0,1]^2\)
\begin{equation}
\begin{split}
& \lim_{v\to 0}
\frac{F\bigl(u+(1-u)v,u\bigr) - F\bigl(u,u+(1-u)v\bigr)}{v} \\
& = (1-u)\,\bigl[F_\xi^{'}(u,u) - F_\eta^{'}(u,u)\bigr].
\end{split}
\end{equation}

We have thus shown that the Duffy-type mapping of the two triangles, followed by their recombination in \eqref{eq:A_I_uv}, transforms the original Cauchy-singular kernel into a regular integrand. The line singularity \(\eta=\xi\) is converted into a simple algebraic factor that is exactly cancelled by the antisymmetric difference of the two evaluations of \(F\). As a consequence, the coincident–element integral \(I\) can be evaluated by tensor–product GL quadrature on the unit square, with the same convergence properties as for a smooth 2D integrand. In particular, no explicit analytical singularity extraction is required, and the argument applies uniformly to any element order.

A subtle issue arises regarding the \emph{discrete} stability of the tensor-product GL quadrature applied to the transformed representation~\eqref{eq:A_I_uv}. Although the integrand in~\eqref{eq:A_I_uv} is regular, one must verify that no numerical instability arises from the clustering of GL nodes near $v = 0$ under the Duffy-type mapping. In particular, each quadrature term contains an effective factor $w_j / v_j$, and it is not \emph{a priori} evident that this ratio remains bounded as $v_j \to 0$.

We now demonstrate that the GL rule is fully stable. Let $\{x_j,w_j\}_{j=1}^n$ denote the GL nodes and weights on $[-1,1]$. Near the endpoint $x=1$, the Mehler--Heine formula \cite{Szego1975} supplies the well--known Bessel correspondence
\begin{equation}
\lim_{n\to\infty} 
P_n\!\left(\cos\frac{z}{n}\right)=\lim_{n\to\infty} 
P_n\!\left(1-\frac{z^2}{2n^2}\right) = J_0(z),
\end{equation}
from which it follows that the nodes near $x=1$ satisfy
\begin{equation}
x_j \;=\; 1 - \frac{j_{0,j}^2}{2n^2} + O(n^{-4}),
\qquad
v_j = \frac{1-x_j}{2}
    = \frac{j_{0,j}^2}{4n^2}+O(n^{-4}),
\label{eq:A_v_asym}
\end{equation}
where $j_{0,j}$ is the $j$-th positive zero of the Bessel function $J_0$. Although the Mehler--Heine formula is an asymptotic statement as $n\to\infty$, the corresponding approximation of the GL nodes near $x=1$ is already good for moderately small orders (e.g., $n\approx 10$). For the GL weights we use the exact identity
\[
w_j=\frac{2}{(1-x_j^2)\,[P_n'(x_j)]^2},
\]
together with the endpoint derivative asymptotic
\rev{\begin{equation}
P_n'(x_j)
= \frac{n^2}{j_{0,j}}\,J_1(j_{0,j}) + O(n),
\end{equation}}
again derived from the Mehler--Heine expansion and its derivative. Substituting these relations yields
\begin{equation}
w_j
= \frac{2}{n^{2}J_1^2(j_{0,j})} + O(n^{-4}).
\label{eq:A_W_asym}
\end{equation}
Combining~\eqref{eq:A_v_asym} and~\eqref{eq:A_W_asym}, the effective factor in \eqref{eq:A_I_uv} becomes
\begin{equation}
\frac{w_j}{v_j}
= \frac{8}{j_{0,j}^{\,2}J_1^2(j_{0,j})} + O(n^{-2}),
\label{eq:A_ratio}
\end{equation}
a finite positive constant depending only on the Bessel zeros. \rev{Although floating-point near-cancellation inevitably occurs in the numerator of \eqref{eq:A_I_uv} as $v_j \to 0$, introducing a local round-off error of $\mathcal{O}(\epsilon_{\mathrm{M}})$, its apparent amplification by the denominator $1/v_j$ is neutralized in the quadrature sum. Specifically, the discrete error contribution from such nodes scales as $w_j \cdot \mathcal{O}(\epsilon_{\mathrm{M}}/v_j)$. Since the ratio $w_j/v_j$ is asymptotically bounded by $\mathcal{O}(1)$ according to \eqref{eq:A_ratio}, the resulting numerical noise is confined to the machine precision $\mathcal{O}(\epsilon_{\mathrm{M}})$.} We conclude that the tensor-product GL quadrature is \emph{asymptotically stable} under the endpoint clustering, fully consistent with the analytical regularity established in \eqref{eq:A_I_uv}.

To complement the analytical regularity and stability proofs, we present a numerical experiment demonstrating that the symmetric Duffy mapping, combined with tensor--product GL rules, yields machine--precision accuracy for Cauchy integrals.
Consider the principal–value integral
\rev{
\begin{equation}
I =
\operatorname{p.v.}
\int_0^1\!\!\int_0^1
\frac{1+\eta}{(1+\xi)(\eta-\xi)}
\,d\xi\,d\eta .
\label{eq:A_rational_test}
\end{equation}
The exact value is
\begin{equation}
I_{\mathrm{exact}} =\ln 2.
\end{equation}}

The integral is evaluated using the coincident–element transformation and standard GL rules in $(u,v)$ coordinates. Table~\ref{tab:appendix_cauchy_test} reports the convergence of the tensor-product Gauss--Legendre rule.
\begin{table}[H]
\centering
{\revcolor
\caption{Numerical evaluation of the Cauchy integral~\eqref{eq:A_rational_test} using tensor-product Gauss--Legendre quadrature after the symmetric split-domain transformation.}
\label{tab:appendix_cauchy_test}
\small
\setlength{\tabcolsep}{8pt}
\renewcommand{\arraystretch}{1.15}
\begin{tabular}{@{} l c c @{}}
\toprule\addlinespace[.25ex]\toprule
\textbf{\(n_q\)} &
\textbf{\(I_{\mathrm{num}}\)} &
\textbf{Error \(\lvert I_{\mathrm{num}}-I_{\mathrm{exact}}\rvert\)} \\
\midrule
2  & \num{0.691082113762525} & \num{2.07e-3}  \\
3  & \num{0.693087878221591} & \num{5.93e-5}  \\
4  & \num{0.693145462290331} & \num{1.72e-6}  \\
5  & \num{0.693147130516645} & \num{5.00e-8}  \\
6  & \num{0.693147179097880} & \num{1.46e-9}  \\
8  & \num{0.693147180558692} & \num{1.25e-12} \\
10 & \num{0.693147180559944} & \num{8.9e-16} \\
12 & \num{0.693147180559945} & \num{1.1e-16} \\
16 & \num{0.693147180559945} & \num{1.1e-16} \\
\bottomrule\addlinespace[.25ex]\bottomrule
\end{tabular}
}
\end{table}

Together with the Mehler--Heine analysis, this example provides strong numerical evidence that the proposed transformation yields
a uniformly stable and high--accuracy discretization of Cauchy-singular integrals for Galerkin BEM.

\section*{Appendix B. Expressions of source field of cylindrical coil}
\setcounter{equation}{0}
\renewcommand{\theequation}{B\arabic{equation}}

\rev{The formulae collected in this appendix are based on the Bessel--Laplace integral identities and elliptic-integral representations of Conway~\cite{Conway2001}, but they are specialized and reorganized here for the cylindrical excitation coil used in the present BEM model. In particular, the normal derivative \(\partial_n A_\varphi^{(e)}\) required by Auld's formula, and the limiting expressions on the symmetry axis are included to make the source term in the impedance calculation directly reproducible.}

Consider a cylindrical coil centered at $(0, 0, z_0)$ in cylindrical coordinates $(r, \phi, z)$, with inner and outer radii $r_1$, $r_2$, axial length $2h$, and $N$ turns. The coil carries an alternating current with amplitude $I$ and angular frequency $\omega$. The electromagnetic field generated by such a coil can be expressed in terms of elliptic integrals~\cite{Conway2001}. Due to the axial symmetry of the configuration, only the $\varphi$-component of the vector potential is non-zero
\begin{equation}
\begin{split}
& A_{\varphi}^{(e)}(r,\zeta)=\frac{\mu_0NI}{4h(r_2-r_1)}
\\
& \times\begin{cases}
\int_{r_1}^{r_2} \rho[f_1(\rho,r,-h-\zeta)-f_1(\rho,r,h-\zeta)] d\rho & \zeta \leq -h\\[1.3ex]
\int_{r_1}^{r_2} \rho[f_1(\rho,r,\zeta-h)-f_1(\rho,r,\zeta+h)] d\rho & \zeta \geq h\\[1.3ex]
g_1(r)-\int_{r_1}^{r_2} \rho[f_1(\rho,r,\zeta+h)+f_1(\rho,r,h-\zeta)] d\rho & -h <\zeta< h
\end{cases}
\end{split}
\end{equation}
where \(\zeta=z-z_0\), and
\begin{equation}
\begin{split}
& f_1(\rho,r,\zeta)=\int_{0}^{\infty}J_1(\lambda r)J_1(\lambda \rho)\frac{e^{-\zeta \lambda}}{\lambda}d\lambda \\
&=\frac{\zeta}{\pi\sqrt{\rho r}}\left[ \frac{\mathbf{E}(\alpha)}{\sqrt{\alpha}}-\frac{(2\rho^2+2r^2+\zeta^2)\sqrt{\alpha}\mathbf{K}(\alpha)}{4\rho r}\right]
\\
& +\frac{1}{2}\begin{cases}
\frac{\rho^2-r^2}{2\rho r}\Lambda_0(\beta,\alpha)+\frac{r}{\rho} & r<\rho \\[1.3ex]
\frac{r^2-\rho^2}{2\rho r}\Lambda_0(\beta,\alpha)+\frac{\rho}{r} & r>\rho
\end{cases}
\end{split}
\label{eq:f1}
\end{equation}
and
\begin{equation}
g_1(r)=\begin{cases}
(r_2-r_1)r & r \leq r_1 \\[1.3ex]
-\frac{2}{3}r^2+r_2r-\frac{r_1^3}{3r} & r_1<r<r_2 \\[1.3ex]
\frac{r_2^3-r_1^3}{3r} & r \geq r_2
\end{cases}
\end{equation}
with
\[
\alpha=\frac{4\rho r}{(\rho+r)^2+\zeta^2}
\]
and
\[
\beta=\arcsin{\frac{\zeta}{\sqrt{(r-\rho)^2+\zeta^2}}}
\]
In equation \eqref{eq:f1}, \(\mathbf{K}(\alpha)\), \(\mathbf{E}(\alpha)\) are the complete elliptic integrals of the first and second kind, \(\Lambda_0\) is the Heuman Lambda function
\[
\Lambda_0(\beta,\alpha)=\frac{2}{\pi}[\mathbf{E}(\alpha)F(\beta,\alpha')+\mathbf{K}(\alpha)E(\beta,\alpha')-\mathbf{K}(\alpha)F(\beta,\alpha')]
\]
where \(\alpha'=1-\alpha\), \(F(\beta,\alpha')\) and \(E(\beta,\alpha')\) are the incomplete elliptic integrals of the first and second kinds, respectively.
The magnetic source field \(B_{z}^{(e)}\) is
\begin{equation}
\begin{split}
& B_{z}^{(e)}(r,\zeta)=\frac{\mu_0 NI}{4h(r_2-r_1)}
\\
& \times\begin{cases}
\int_{r_1}^{r_2} \rho[f_2(\rho,r,-h-\zeta)-f_2(\rho,r,h-\zeta)] d\rho & \zeta \leq -h\\[1.3ex]
\int_{r_1}^{r_2} \rho[f_2(\rho,r,\zeta-h)-f_2(\rho,r,\zeta+h)] d\rho & \zeta \geq h\\[1.3ex]
g_2(r)-\int_{r_1}^{r_2} \rho[f_2(\rho,r,\zeta+h)+f_2(\rho,r,h-\zeta)] d\rho & -h <\zeta< h
\end{cases}
\end{split}
\end{equation}
where
\begin{equation}
\begin{split}
f_2(\rho,r,\zeta)&=\int_{0}^{\infty}J_0(\lambda r)J_1(\lambda \rho)e^{-\zeta \lambda}d\lambda \\
&=\frac{1}{\rho}\cdot\begin{cases}
1-\frac{\zeta\sqrt{\alpha}\mathbf{K}(\alpha)}{2\pi\sqrt{\rho r}}-\frac{\Lambda_0(\beta,\alpha)}{2} & r<\rho \\[1.3ex]
-\frac{\zeta\sqrt{\alpha}\mathbf{K}(\alpha)}{2\pi\sqrt{\rho r}}+\frac{\Lambda_0(\beta,\alpha)}{2} & r>\rho
\end{cases}
\end{split}
\end{equation}
and
\begin{equation}
g_2(r)=\begin{cases}
2(r_2-r_1) & r \leq r_1 \\[1.3ex]
2(r_2-r) & r_1<r<r_2 \\[1.3ex]
0 & r \geq r_2
\end{cases}
\end{equation}
The magnetic source field \(B_{r}^{(e)}\) is
\begin{equation}
B_{r}^{(e)}(r,\zeta)=\frac{\mu_0 NI}{4h(r_2-r_1)} \int_{r_1}^{r_2} \rho \left[f_3(\rho,r,\lvert h-\zeta \rvert)-f_3(\rho,r,\lvert h+\zeta \rvert)\right] d\rho
\end{equation}
where
\begin{equation}
\begin{split}
f_3(\rho,r,\zeta)&=\int_{0}^{\infty}J_1(\lambda r)J_1(\lambda \rho)e^{-\zeta \lambda}d\lambda \\
&=\frac{1}{\pi\sqrt{\alpha \rho r}}\left[(2-\alpha)\mathbf{K}(\alpha)-2\mathbf{E}(\alpha)\right]
\end{split}
\end{equation}
For the Auld's formula, the normal derivative of the vector potential is
\begin{equation}
\partial_{n}A_{\varphi}^{(e)}=n_r\partial_r A_{\varphi}^{(e)}+n_z\partial_zA_{\varphi}^{(e)}
\end{equation}
where
\begin{equation}
\partial_r A_{\varphi}^{(e)}=B_z^{(e)}-A_{\varphi}^{(e)}/r
\end{equation}
and
\begin{equation}
\partial_zA_{\varphi}^{(e)}=-B_r^{(e)}=-\frac{\mu_0 NI}{4h(r_2-r_1)}\int_{r_1}^{r_2}\rho[f_3(\rho,r,|h-\zeta|)-f_3(\rho,r,|h+\zeta|)]d\rho
\end{equation}
When \(r=0\), the above field components reduce to
\begin{equation}
A_{\varphi}^{(e)}(0,\zeta)=0
\end{equation}
\begin{equation}
B_{r}^{(e)}(0,\zeta)=0
\end{equation}
\begin{equation}
\begin{split}
& B_{z}^{(e)}(0,\zeta)=\frac{\mu_0 NI}{4h(r_2-r_1)}
\\
& \times \left[(\zeta+h)\ln{\frac{r_2+\sqrt{r_2^2+(\zeta+h)^2}}{r_1+\sqrt{r_1^2+(\zeta+h)^2}}}-(\zeta-h)\ln{\frac{r_2+\sqrt{r_2^2+(\zeta-h)^2}}{r_1+\sqrt{r_1^2+(\zeta-h)^2}}}\right]
\end{split}
\end{equation}
and
\begin{equation}
\partial_rA_{\varphi}^{(e)}(0,\zeta)=\frac{1}{2}B_{z}^{(e)}(0,\zeta)
\end{equation}

\section*{Appendix C. Axisymmetric kernel of a ring source}
\setcounter{equation}{0}
\renewcommand{\theequation}{C\arabic{equation}}

The axisymmetric Laplace kernel $\mathcal{G}$ is
\begin{equation}
\begin{split}
\mathcal{G}(r,z;r',z')&=\frac{1}{4\pi}\int_{0}^{2\pi}\frac{\cos(\varphi)}{R(\varphi)}d\varphi \\
&=\frac{1}{2}\int_{0}^{\infty}J_1(\lambda r)J_1(\lambda r')e^{-\lambda|z-z'|}d\lambda \\
&= \frac{1}{\pi \sqrt{m r r'}} \left[ \left(1-\frac{m}{2}\right) K(m) - E(m) \right].
\end{split}
\label{eq:axisym_laplace_kernel}
\end{equation}
where
\begin{equation}
R(\varphi)=\sqrt{r^2+r'^2-2rr'\cos(\varphi)+(z-z')^2}
\label{eq:axisym_distance_R}
\end{equation}
and
\begin{equation}
m = \frac{4rr'}{(r+r')^2 + (z-z')^2}.
\label{eq:axisym_modulus_m}
\end{equation}
The components for partial derivatives are
\begin{equation}
\begin{split}
\partial_{r'}\mathcal{G}=\frac{r}{2\pi(mrr')^{3/2}}\left\{m\left[\left(\frac{m}{2}-1\right)K(m)+E(m)\right]\right.\\
\left.+\frac{m_{r'}r'}{m-1}\left[(1-m)K(m)+\left(\frac{m}{2}-1\right)E(m)\right]\right\}
\end{split}
\label{eq:laplace_kernel_drprime}
\end{equation}
and
\begin{equation}
\begin{split}
\partial_{z'}\mathcal{G}=\frac{m_{z'}rr'}{2\pi(m-1)(mrr')^{3/2}}\left[(1-m)K(m)+\left(\frac{m}{2}-1\right)E(m)\right].
\end{split}
\label{eq:laplace_kernel_dzprime}
\end{equation}
where
\begin{equation}
m_{r'}=\frac{4r[r^2-r'^2+(z-z')^2]}{[(r+r')^2+(z-z')^2]^2}
\label{eq:m_rprime}
\end{equation}
and
\begin{equation}
m_{z'}=\frac{8rr'(z-z')}{[(r+r')^2+(z-z')^2]^2}.
\label{eq:m_zprime}
\end{equation}
When $r=0$, we have $J_1(\lambda r)=0$. Hence, from the Bessel integral expression it is easy to conclude that 
\[
\mathcal{G}(0,z;r',z')=0,\quad
\partial_{r'}\mathcal{G}(0,z;r',z')=0,\quad
\partial_{z'}\mathcal{G}(0,z;r',z')=0.
\]
When $r'=0$, we have $\mathcal{G}(r,z;0,z')=0$, $\partial_{z'}\mathcal{G}(r,z;0,z')=0$, and
\begin{equation}
\partial_{r'}\mathcal{G}(r,z;0,z')=\frac{r}{4\chi^{3}}
\label{eq:laplace_axis_rprime_limit}
\end{equation}
where
\begin{equation}
\chi=\sqrt{r^2+(z-z')^2}.
\label{eq:axisym_chi}
\end{equation}

The axisymmetric Helmholtz kernel $\mathcal{G}_k$ is
\begin{equation}
\begin{split}
\mathcal{G}_k&=\frac{1}{4\pi}\int_{0}^{2\pi}\frac{e^{-ikR(\varphi)}}{R(\varphi)}\cos(\varphi)d\varphi \\
&=\frac{1}{2}\int_{0}^{\infty}J_1(\lambda r)J_1(\lambda r')\frac{\lambda}{\sqrt{\lambda^2-k^2}}e^{-\sqrt{\lambda^2-k^2}|z-z'|}d\lambda .
\end{split}
\label{eq:axisym_helmholtz_kernel}
\end{equation}
The Bessel integral expression \eqref{eq:axisym_helmholtz_kernel} is derived from \cite{Sommerfeld1978} and the Graf's addition theorem \cite{Watson1944}. In practice, we use singularity extraction
\begin{equation}
\mathcal{G}_k=\frac{1}{4\pi}\int_{0}^{2\pi}\frac{e^{-ikR(\varphi)}-1}{R(\varphi)}\cos(\varphi)d\varphi+\mathcal{G}.
\label{eq:helmholtz_extraction}
\end{equation}
The integral corresponding to the first term on the right-hand side of \eqref{eq:helmholtz_extraction} can be evaluated accurately using GL quadrature. The second term contains the Laplace kernel, whose singularity requires the regularization technique in \textsection~3. An analogous approach applies to the partial derivatives of the Helmholtz kernel:
\begin{equation}
\begin{split}
\partial_{r'}\mathcal{G}_k&=\frac{1}{4\pi}\int_{0}^{2\pi}\frac{(r'-r\cos(\varphi))\cos(\varphi)\left[-(1+ikR(\varphi))e^{-ikR(\varphi)}+1\right]}{R^3(\varphi)}d\varphi\\
& +\partial_{r'}\mathcal{G}
\end{split}
\label{eq:helmholtz_drprime_extraction}
\end{equation}
and
\begin{equation}
\begin{split}
\partial_{z'}\mathcal{G}_k&=\frac{1}{4\pi}\int_{0}^{2\pi}\frac{(z-z')\cos(\varphi)\left[(1+ikR(\varphi))e^{-ikR(\varphi)}-1\right]}{R^3(\varphi)}d\varphi\\
& +\partial_{z'}\mathcal{G}.
\end{split}
\label{eq:helmholtz_dzprime_extraction}
\end{equation}
The static terms on the right-hand sides of \eqref{eq:helmholtz_drprime_extraction} and \eqref{eq:helmholtz_dzprime_extraction}, defined in \eqref{eq:laplace_kernel_drprime} and \eqref{eq:laplace_kernel_dzprime}, can be treated using the same regularization procedure. In the limit $r \to 0$, the asymptotic behavior follows directly from the Bessel function representation in \eqref{eq:axisym_helmholtz_kernel}:
\[
\mathcal{G}_k(0,z;r',z')=0,\quad
\partial_{r'}\mathcal{G}_k(0,z;r',z')=0,\quad
\partial_{z'}\mathcal{G}_k(0,z;r',z')=0.
\]
When $r'=0$, we have $\mathcal{G}_k(r,z;0,z')=0$, $\partial_{z'}\mathcal{G}_k(r,z;0,z')=0$, and
\begin{equation}
\begin{split}
\partial_{r'}\mathcal{G}_k(r,z;0,z')&=\frac{1}{4}\int_{0}^{\infty}J_1(\lambda r)\frac{\lambda^2}{\sqrt{\lambda^2-k^2}}e^{-\sqrt{\lambda^2-k^2}|z-z'|}d\lambda\\
&=\frac{r}{4\chi^3}(1+ik\chi)e^{-ik\chi}.
\end{split}
\label{eq:helmholtz_axis_rprime_limit}
\end{equation}

\bibliographystyle{elsarticle-num}
\bibliography{Refs}

\end{document}